\documentclass{article}
\usepackage[pdftex]{hyperref}
\usepackage{amsmath,amsfonts,amsthm,amssymb,graphicx,tikz,tikz-cd}
\usepackage[all,2cell,ps]{xy}

\theoremstyle{plain}
\newtheorem{thm}{Theorem}[section]
\newtheorem{lem}[thm]{Lemma}
\newtheorem{prop}[thm]{Proposition}
\newtheorem{cor}[thm]{Corollary}

\theoremstyle{definition}
\newtheorem{defn}[thm]{Definition}

\newcommand{\Z}{\mathbb Z}

\newcommand{\Q}{\mathbb Q}

\newcommand{\C}{\mathbb C}
\newcommand{\pp}{\mathbb{P}}

\newcommand{\lam}{\lambda}

\newenvironment{pf}{\begin{proof}}{\end{proof}}

\usetikzlibrary{arrows}

\title{The Toledo invariant, and Seshadri constants of fake projective planes}
\author{Luca F.\ Di Cerbo\footnote{This material is based upon work supported by a Grant of the Max Planck Society: ``Complex Hyperbolic Geometry and Toroidal Compactifications'' and a grant associated to the S.S. Chern position at ICTP.} \\ \small{Abdus Salam International Centre for Theoretical Physics - ICTP} \\ 
\small{\textsf{ldicerbo@ictp.it}}}

\begin{document}

\maketitle

\begin{abstract}
The purpose of this paper is to explicitly compute the Seshadri constants of all ample line bundles on fake projective planes. The proof relies on the theory of the Toledo invariant, and more precisely on its characterization of $\C$-Fuchsian curves in complex hyperbolic spaces.
\end{abstract}

\section{Introduction and Preliminaries}\label{sec:Intro}

Seshadri constants measure the local positivity of an ample line bundle on a projective variety. They were introduced by Demailly in his study of Fujita's conjecture through the theory of positive currents and singular Hermitian metrics, see \cite{Dem90}. Given a nef line bundle $L$ and a point $x\in X$, let us define
\begin{align}\notag
\epsilon(L, x)=\inf_{C\supset x}\frac{L\cdot C}{\text{mult}_{x}(C)},
\end{align} 
to be the Seshadri constant of $L$ at $x$, where $C\subset X$ is an irreducible holomorphic curve and $\text{mult}_{x}(C)$ is the multiplicity of such a curve at $x$. The global Seshadri constant of $L$ is then given by 
\begin{align}\notag
\epsilon(L)=\inf_{x\in X}\epsilon(L, x),
\end{align}
where by the Seshadri criterion for ampleness, there is a $\epsilon>0$ such that $\epsilon(L)>\epsilon$ if and only if $L$ is ample, see for example Theorem 1.4.13 in \cite{Laz1}. For more on this circle of ideas we refer to the survey \cite{Di Rocco}.

Despite being very easy to define, the Seshadri constants are in practice hard to compute or even estimate. In recent years, there has been a growing interest in the study of Seshadri constants on locally symmetric spaces via geometric and analytic techniques. For example, Hwang and To in \cite{Hwang-To} were able to estimate from below the Seshadri constants, at any given point, of the canonical line bundle of a compact complex hyperbolic space, in terms of the radius of largest geodesic ball centered at that point with respect to the locally symmetric Bergman metric. For many more results in this direction, we refer to \cite{Di Rocco} and the bibliography therein.

The purpose of this paper is to compute the Seshadri constants of all ample line bundles on fake projective planes. 
Let us recall their definition:
\begin{defn}
A fake projective plane is a surface of general type $X$ with $c_{2}=3$ and $p_{g}=H^{0}(X; K_{X})=0$.
\end{defn}

These surfaces are commonly referred as ``fake projective planes'', since it can be easily shown they must have the same rational homology groups as $\pp^{2}$. Fake projective planes are particular compact complex hyperbolic surfaces, see Section \ref{FPP} for more details. It turns out that, given any ample line bundle $L$ on a fake projective plane $X$, the pointwise Seshadri constants are independent of the point, so that $\epsilon(L, x)=\epsilon(L)$ for any $x\in X$. This fact is particularly interesting as fake projective planes are \emph{not} homogeneous, even if they are clearly locally symmetric. In fact, roughly half of all fake projective planes indeed have trivial automorphism group. Finally, we also explicitly compute $\epsilon(L)$ for all ample line bundles. For the precise numerical values, we refer to Theorem \ref{main} in Section \ref{final}.

The paper is organized as follows. Section \ref{FPP} starts with the definition of a fake projective plane and collects some of the basic properties of these surfaces. More precisely, we recall some important features of their fundamental group and homology which follow from the Prasad-Yeung and Cartwright-Steger classification. In Section \ref{Fuchsian}, we recall the basic theory of the Toledo invariant. We use this invariant to study curves on fake projective planes, and in particular to understand their singularities. Then, in Section \ref{final} we prove Theorem \ref{main} which computes exactly the Seshadri constants of a fake projective plane.\\

\noindent\textbf{Acknowledgments}. I would like to thank the Max Planck Institute for Mathematics for the great working environment while this project was conceived and completed. I would also like to thank Matthew Stover for introducing me to the beautiful ideas of Toledo.

\subsection{Fake Projective Planes}\label{FPP}

In this section, we recall some of the general properties of the so-called fake projective planes. For the basic complex surface theory and complex hyperbolic geometry, we refer to the books \cite{Bar}, \cite{Bea} and \cite{Gold Book}.

The formula of Noether for the holomorphic Euler characteristic of a smooth surface gives
\[
\chi_{\mathcal{O}}=1-q+p_{g}=\frac{c^{2}_{1}+c_{2}}{12},
\]
where $q=h^{1,0}$ is the irregularity of the surface. Since $\chi_{\mathcal{O}}>0$ for any surface of general type, we conclude that a fake projective plane must have vanishing first Betti number, i.e., $q=0$ which is equivalent to $H_{1}(X; \Q)=0$. Thus, the fact $c_{2}=3$ implies that $h^{1, 1}=1$, so that the Picard number of any such space is always one. Moreover, let us observe that any fake projective plane achieves the equality sign in the so-called Bogomolov-Miyaoka-Yau inequality
\[
c^{2}_{1}\leq 3c_{2},
\]
so that, because of Yau's solution of Calabi conjecture \cite{Yau}, it is a compact complex hyperbolic $2$-manifold.
Thus, any fake projective plane $X$ is given as the quotient of the complex hyperbolic $2$-space $\mathcal{H}^{2}_{\C}$ by a torsion free co-compact lattice $\Gamma\in \text{PU}(2, 1)$. Furthermore, Klingler \cite{Klingler} and Yeung \cite{Yeung} proved that any such $\Gamma$ is always an arithmetic subgroup of $\text{PU}(2, 1)$. This important result is crucial for the classification of Prasad-Yeung \cite{Prasad} and Cartwright-Steger \cite{Steger} of all fake projective planes. More precisely, there are 100 isomorphism classes of fake projective planes and 50 explicit subgroups in $\text{PU}(2, 1)$ which arise as fundamental groups of fake projective planes. For any fake projective plane $X$ with $\pi_{1}(X)=\Gamma$, we then have that
\[
H_{1}(X; \Z)=\Gamma/[\Gamma, \Gamma]
\]
is a finite abelian group which is never zero, see again \cite{Prasad} and also the addendum \cite{Prasad2}. Moreover, since $H^{1}(X; \mathcal{O})=0$ we have that $Pic(X)=H^{2}(X; \Z)$, and by the universal coefficient theorem $\text{Tor}(H^{2}(X; \Z))=H_{1}(X; \Z)\neq 0$. Thus, on any fake projective plane we always have torsion line bundles. Nevertheless,  by Poincar\'e duality $Pic(X)/\text{Tor}(H^{2}(X;\Z))$ is one dimensional unimodular lattice and we therefore have the existence of an ample line bundle $L$ such that $c_{1}(L)^{2}=1$, which generates the torsion free part of the Picard group. This fact motivates the following definition.

\begin{defn}
For a fake projective plane X, we denote by $L_{1}$ any ample generator of the torsion free part of $Pic(X)$. 
\end{defn}

Let us remark that the choice of $L_{1}$ is not unique. In fact, as we discussed above the different choices are parametrized by $H_{1}(X; \Z)$. Nevertheless, this lack of non-uniqueness does not play an important role for us. This is because the Seshadri constants of an ample line bundle depend only on its numerical equivalence class. In other words, if we denote numerical equivalence by $\equiv$, given any two line bundles $L$ and $L'$ such that $L\equiv L'$, we then have $\epsilon(L, x)=\epsilon(L', x)$ for any point $x$. Next, given any ample line bundle $L$ with self-intersection $k^{2}$ let us observe that
\[
L\equiv kL_{1},
\]
where $L_{1}$ is any ample generator of the torsion free part of $Pic(X)$. Thus, for our purposes, it is natural to give the following definition.

\begin{defn}
For a fake projective plane X, we denote by $L_{k}$ any ample line bundle $L$ such that $c^{2}_{1}(L)=k^{2}$. 
\end{defn}

\subsection{Toledo Invariant and $\C$-Fuchsian curves in ball quotients}\label{Fuchsian}

Let $C$ be a closed hyperbolic Riemann surface and let
\[
\rho: \pi_{1}(C)\rightarrow \text{PU}(2, 1)
\] 
be a representation of its fundamental group into the group of holomorphic isometries of the complex hyperbolic $2$-space. The representation $\rho$ determines a flat $\mathcal{H}^{2}_{\C}$ bundle over $C$. The complex hyperbolic $2$-space can be thought as the unit ball in $\C^{2}$ which is contractible. Thus, this bundle has a section and then there exists a $\rho$-equivariant map
\[
f: \mathcal{H}\rightarrow \mathcal{H}^{2}_{\C}
\] 
where $\mathcal{H}$ is the universal cover of $C$. Let us denote by $\omega_{1}$ the Bergman metric on $\mathcal{H}^{2}_{\C}$ with constant negative holomorphic sectional curvature normalized to be $-1$. Since $f^{*}\omega_{1}$ is invariant under the action of $\pi_{1}(C)$ on $\mathcal{H}$, we define the number
\[
T(\rho)=\frac{1}{2\pi}\int_{C}f^{*}\omega_{1}
\]
to be the Toledo invariant of the representation $\rho$. Note that the Toledo invariant is independent of the map $f$.

\begin{thm}[D. Toledo, \cite{Toledo}]\label{Toledo}
Let $\rho: \pi_{1}(C)\rightarrow \textrm{PU}(2, 1)$ be a representation with $C$ is a closed hyperbolic Riemann surface. We then have 
\[
|T(\rho)|\leq 2g(C)-2
\]
with equality if and only if $\rho$ is $\C$-Fuchsian.
\end{thm}

Recall that $\rho: \pi(C)\rightarrow \text{PU}(2, 1)$ is said to be $\C$-Fuchsian if the $\rho$-equivariant map $f$ is a holomorphic totally geodesic embedding of the complex hyperbolic $1$-space $\mathcal{H}^{1}_{\C}$ into $\mathcal{H}^{2}_{\C}$.\\

For more on this circle of ideas, we refer the interested reader to the paper of Goldman-Kapovich-Leeb \cite{Goldman} and to the recent survey of Stover \cite{Stover} and the bibliography therein.

\section{Curves on ball quotient surfaces and fake projective planes}\label{ssec:Balll}

In this section, we derive a bound on the genus of the normalization of a holomorphic curve in a ball quotient surfaces. This bound follows from the bound on the Toledo invariant presented in Section \ref{Fuchsian}. Moreover, we use this bound to study the geometry of curves on fake projective planes.

To this aim, let $(X, \omega_{1})$ be smooth complex hyperbolic surface equipped with its locally symmetric Bergman metric $\omega_{1}$. Moreover, let us normalize the holomorphic sectional curvature of $\omega_{1}$ to be $-1$. Given a holomorphic curve $C\subset X$, let us consider the normalization map $i: \overline{C}\rightarrow X$. Recall that this map is simply the resolution of the singularities of $C$, so that $\overline{C}$ is a smooth curve and the map $i$ is generically the identity. Moreover, the holomorphic map $i$ is the identity if and only if the curve $C$ is smooth. Since $X$ is hyperbolic, we must have $g(\overline{C})\geq 2$ for any curve $C\subset X$. If otherwise, we would have a non trivial holomorphic map $h: \C\rightarrow X$. The existence of such map contradicts the hyperbolicity of $X$.

Thus given any $C\subset X$, let us consider the normalization map $i:\overline{C}\rightarrow X$. Since $g(\overline{C})\geq 2$, we can consider the Toledo invariant associated to the normalization map:
\[
T(i)=\frac{1}{2\pi}\int_{\overline{C}}i^{*}\omega_{1}.
\]
In this case it is clear that $T(i)>0$, so that by Theorem \ref{Toledo} we have
\[
0<T(i)\leq 2g(\overline{C})-2
\]
with equality if and only if $i$ lifts to a holomorphic totally geodesic embedding of $\mathcal{H}^{1}_{\C}$ into $\mathcal{H}^{2}_{\C}$. In other words, the equality sign is achieved if and only if $C$ is a totally geodesic immersed curve in $X$. 

Next, let us observe the following. Since $(X, \omega_{1})$ is Einstein and the holomorphic sectional curvature is normalized to be $-1$, we have 
\[
c_{1}(K_{X})=\frac{3}{4\pi}\omega_{1},
\]
where by $c_{1}(K_{X})$ we denote the first Chern class of the canonical line bundle $K_{X}$.
Thus, given any curve $C\subset X$ and denoting by $C^{*}$ its smooth locus, we compute
\[
K_{X}\cdot C=\int_{C^{*}}\frac{3}{4\pi}\omega_{1}=\frac{3}{4\pi}\int_{\overline{C}}i^{*}\omega_{1}=\frac{3}{2}T(i).
\]
We can then use the bound on the Toledo invariant of the map $i:\overline{C}\rightarrow X$ to derive the following upper bound for the intersection of $K_{X}$ with $C$.

\begin{prop}\label{toledo}
Let $X$ be a complex hyperbolic surface. Given a reduced irreducible curve $C\in X$, let us denote by $\overline{C}$ its normalization. We then have
\[
K_{X}\cdot C\leq 3(g(\overline{C})-1)
\]
with equality if and only if $C$ is an immersed totally geodesic curve.
\end{prop}

Next, we want to refine this proposition when $X$ is a fake projective plane. Let us start with a well known lemma. We sketch the proof for the convenience of the reader.

\begin{lem}\label{nogeodesic}
There are no immersed totally geodesic curves in a fake projective plane.
\end{lem}

\begin{pf}
The main theorem in \cite{Prasad} combined with \cite{Steger} tells us that any fake projective plane is an arithmetic ball quotient of the second type. These particular ball quotients are known not to carry any immersed totally geodesic curve, see for example page 901 in \cite{Moller-Toledo}. The proof is then complete.
\end{pf}

Thus, let $C$ be any reduced irreducible curve in a fake projective plane $X$. Any such $C$ is numerically equivalent to $kL_{1}$ for some integer $k\geq 1$, where $L_{1}$ is any ample generator of the torsion free part of $Pic(X)$, see Section \ref{FPP}. In particular, we always have that $C\equiv L_{k}$ for some $k\geq 1$. 

\begin{prop}\label{key}
Let $C$ be a reduced irreducible curve in a fake projective plane $X$ numerically equivalent to $L_{k}$ for some $k\geq 1$. Let $\overline{C}$ be its normalization. We then have $g(\overline{C})>1+k$.
\end{prop}
\begin{pf}
Given any fake projective plane $X$, recall that $c^{2}_{1}(X)=K^{2}_{X}=9$. Thus, $K_{X}$ is numerically equivalent to $3L_{1}$. Since $C$ is assumed to be numerically equivalent to $L_{k}$,
by using Proposition \ref{Toledo} we compute that
\[
K_{X}\cdot C=3k\leq3(g(\overline{C})-1)\quad \Rightarrow\quad g(\overline{C})\geq 1+k.
\]
Now the equality sign can be achieved if and only if $C$ is an immersed totally geodesic curve. By Lemma \ref{nogeodesic}, we know that such curves cannot exist on a fake projective plane. We then conclude that
\[
g(\overline{C})>1+k,
\]
as claimed.
\end{pf}

It is interesting to explicitly state a corollary of Proposition \ref{key} for curves of fake projective planes numerically equivalent to $L_{1}$.

\begin{cor}\label{smooth}
Let $C$ be a curve on a fake projective plane $X$ numerically equivalent to $L_{1}$. We then have that $C$ is necessarily a smooth genus three curve.
\end{cor} 

\begin{pf}
By Proposition \ref{key}, we know that $C$ has to be smooth. Finally, we compute that 
\[
g(C)=p_{a}(C)=1+\frac{K_{X}\cdot C+C^{2}}{2}=3
\]
where $p_{a}(C)$ is the arithmetic genus of $C$.
\end{pf}

Let us point out that curves of low degree on fake projective planes are quite mysterious objects. In particular, it seems currently unknown weather or not there exist curves linearly equivalent to $L_{i}$, for $i=1,2$, on any of the 50 fake projective planes.  In any case, if curves numerically equivalent to $L_{1}$ exist on some fake projective plane, they necessarily must be smooth.\\

Next, we want to apply a similar argument derive an upper bound for the multiplicity of a singular point of a reduced, irreducible curve $C\subset X$, which satisfies $C\equiv L_{k}$ for some $k\geq 2$. This upper bound is the key fact for the computation of the Seshadri constants of any line bundle $L_{k}$, on any fake projective plane $X$.

\begin{prop}\label{bound}
Let $C$ be a reduced, irreducible singular curve in a fake projective plane $X$. Let $C$ be numerically equivalent to $L_{k}$ for some $k\geq 2$. For any singular point $p\in C$, we have $2\leq m_{p}\leq k$, where $m_{p}$ denotes the multiplicity of $p$.
\end{prop}

\begin{pf}
Given a singular point $p\in C$, let us define the integer $m_{p}(C)=m_{p}$ to be its multiplicity. Thus $m_{p}\geq 2$, and let us consider the normalization map $i: \overline{C}\rightarrow C$. 
The local genus drop at the point $p$ is defined to be
\[
\delta_{p}=\text{dim}_{\C}(i_{*}\mathcal{O}_{\overline{C}}/\mathcal{O}_{C})_{p}
\]
where by construction we have $g(\overline{C})=p_{a}(C)-\sum_{i}\delta_{p_{i}}$, with the sum taken over all the singular points, say $\{p_{i}\}$, of $C$.  For any singular point $p\in C$ we have that $\delta_{p}\geq m_{p}(m_{p}-1)/2$. 
Thus, we conclude
\[
p_{a}(C)-\frac{m_{p}(m_{p}-1)}{2}\geq p_{a}(C)-\sum_{i}\frac{m_{p_{i}}(m_{p_{i}}-1)}{2}\geq g(\overline{C}).
\]
Since $C$ is numerically equivalent to $L_{k}$, we have that
\[
p_{a}(C)=1+\frac{3k+k^{2}}{2}
\]
so that using the bound given in Proposition \ref{key}, we obtain
\[
3k+k^{2}-m_{p}(m_{p}-1)>2k\quad \Rightarrow \quad m_{p}^{2}-m_{p}-k-k^{2}<0.
\]
In other words, for any singular point $p\in C$, we proved that 
\[
2\leq m_{p}< 1+k,
\]
which concludes the proof.
\end{pf}

\section{Seshadri constants of fake projective planes}\label{final}

In this section, we finally compute the Seshadri of all line bundles $L_{k}$ on any of the fake projective planes.  Let us start with a proposition which is a consequence of a general result of Ein-Lazarsfeld \cite{Ein} for Seshadri constants on smooth surfaces. The main result of this section, Theorem \ref{main} below, is sharp generalization of this proposition.

\begin{prop}\label{Laz}
Let $X$ be a fake projective plane. Let $L_{k}$ be an ample line bundle with self-intersection $k^{2}$. We then have $\epsilon(L_{k}, x)=k$ for all except at most countably many points $x\in X$.
\end{prop}
\begin{pf}
Given $L_{k}$ recall that $L_{k}\equiv kL_{1}$. Now $L_{1}$ is an ample line bundle with $L_{1}^{2}=1$. By a result of Ein-Lazarsfeld
\cite{Ein}, we have that $\epsilon(L_{1}; x)=1$ for all except possibly countably many points in $X$. Since $\epsilon(L_{k}; x)=k\epsilon(L_{1}; x)$, the proof is complete.
\end{pf}

We can now prove the main theorem.

\begin{thm}\label{main}
Let $X$ be a fake projective plane. Let $L_{k}$ be an ample line bundle with self-intersection $k^{2}$. Given any point $x\in X$, we have $\epsilon(L_{k}, x)=\epsilon(L_{k})=k$.
\end{thm}
\begin{pf}
Let $C$ be a curve in $X$ numerically equivalent to $L_{1}$. By Corollary \ref{smooth}, any such curve is smooth so that
\[
\frac{L_{k}\cdot C}{\text{mult}_{p}(C)}=k
\]
for any point $p\in C$. Next, let $C$ be a curve in $X$ numerically equivalent $L_{l}$ for some $l\geq 2$. By Proposition \ref{bound}, if $p\in C$ is a point with multiplicity $m_{p}\geq 2$ we have the bound $m_{p}\leq l$. Thus, we then compute
\[
\frac{L_{k}\cdot C}{\text{mult}_{p}(C)}\geq \frac{l\cdot k}{l}=k.
\] 
In conclusion, for any point $x\in X$, we have the inequality $\epsilon(L_{k}; x)\geq k$. In order to finish the proof, we need to show that $\epsilon(L_{k}; x)\leq k$ for any $x\in X$. To this aim, let $f:Y\rightarrow X$ be the blow up map at $x\in X$ and let us denote by $E$ the exceptional divisor in $Y$. Let us observe that the Seshadri constant of $L_{k}$ at $x$ can be equivalently defined as follows \cite{Laz1}: 
\[
\epsilon(L_{k}, x)=\text{sup}\{\lam>0  :  f^{*}L_{k}-\lam E\quad \text{is nef on Y} \}.
\]
For any $x\in X$ we then have
\[
\epsilon(L_{k}, x)\leq \sqrt{L^{2}_{k}}=k.
\]
The proof is complete.
\end{pf}

In particular, we can explicitly compute the Seshadri constants of the canonical line bundle of a fake projective plane.

\begin{cor}\label{canonical}
Let $X$ be a fake projective plane. Given any point $x\in X$, we have $\epsilon(K_{X}, x)=\epsilon(K_{X})=3$.
\end{cor}
\begin{pf}
By Theorem \ref{main}, it suffices to observe that $K_{X}\equiv L_{3}$.
\end{pf}

It would be interesting to compare the estimates given by  Hwang and To in \cite{Hwang-To} for $\epsilon(K_{X})$, with the exact result given in Corollary \ref{canonical}. This comparison would require the computation of the injectivity radii of fake projective planes, which seems an interesting problem on its own.\\

Let us conclude this section by proving a result on the $3$-canonical map of a fake projective plane. More precisely, we show that the map associated to the linear system $|3K_{X}|$ gives an embedding. The result is a direct consequence Corollary \ref{canonical}.

\begin{cor}\label{canonical2}
Let $X$ be a fake projective plane. The $3$-canonical map 
\[
\varphi_{|3K_{X}|}: X\rightarrow \pp^{27}
\]
is an embedding.
\end{cor}
\begin{pf}
By Corollary \ref{canonical}, we have $\epsilon(2K_{X})=6>2\cdot\text{dim}_{\C}(X)=4$. Thus, using Proposition 6.8 in \cite{Dem90} we know that the map associated to the linear system $|K_{X}+2K_{X}|$ gives an embedding into some projective space. It remains to compute the dimension of $H^{0}(X; 3K_{X})$. Using the vanishing theorem of Kodaira
and the Riemann-Roch formula, we obtain $h^{0}(X; 3K_{X})=28$.
\end{pf}

Let us remark that Corollary \ref{canonical2} can alternatively be proved by using Reider's theorem, see for example page 176 in \cite{Bar}. This fact was already observed in Section 10 of \cite{Prasad}. Nevertheless, it seems of interest to present a self-contained Seshadri constants proof based on the exact result given in Corollary \ref{canonical}.



\begin{thebibliography}{ELMNPM}

\bibitem[BHPV04]{Bar} W. P. Barth, K. Hulek, C. A. Peters, A. Van de Ven, Compact complex surfaces. Second edition. Ergebnisse der Mathematik und ihrer Grenzgebiete. 3. Folge. A Series of Modern Surveys in Mathematics [Results in Mathematics and Related Areas. 3rd Series. A Series of Modern Surveys in Mathematics], 4. \textit{Springer-Verlag, Berlin}, 2004.

\bibitem[Bauer et al.09]{Di Rocco} T. Bauer, S. Di Rocco, B. Harbourne, M. Kapustka, A. Knutsen, W. Syzdek, T. Szemberg, A primer on Seshadri constants. \textit{Interactions of classical and numerical algebraic geometry}, 33-70, Contemp. Math., 496, \textit{Amer. Math. Soc., Providence, RI,} 2009.

\bibitem[Bea96]{Bea} A. Beauville, Complex Algebraic Surfaces. Second Edition,
London Mathematical Society Students Texts, 34. \textit{Cambridge University Press, Cambridge}, 1996.

\bibitem[CS10]{Steger} D. Cartwright, T. Steger, Enumeration of the 50 fake projective planes.
\textit{C. R. Acad. Sci. Paris, Ser. I} \textbf{348} (2010), 11-13.

\bibitem[Dem90]{Dem90} J.- P. Demailly, Singular Hermitian metrics on positive line bundles.
\textit{Complex Algebraic Varieties (Bayreuth, 1990)}, 87-104,  Lecture Notes in Math., 1507, \textit{Springer, Berlin,} 
1992. 

\bibitem[EL93]{Ein} L. Ein, R. Lazarsfeld,  Seshadri constants on smooth surfaces.
Journ\'ees de G\'eom\'etrie Alg\'ebrique d'Orsay (Orsay, 1992).
\textit{Ast\'erisque} \textit{No. 218} (1993), 177-186.

\bibitem[Gol99]{Gold Book} William R. Goldman,
Complex hyperbolic geometry. Oxford Mathematical Monographs. Oxford Science Publications. \textit{The Claredon Press, Oxford University Press, Oxford University Press, New York,}, 1999.

\bibitem[GKL01]{Goldman} W. Goldman, M. Kapovich, B. Leeb, Complex hyperbolic manifolds homotopy equivalent to a Riemann surface. \textit{Comm. Anal. Geom.} \textbf{9} (2001), no.1, 61-95.

\bibitem[HT99]{Hwang-To} J.- M. Hwang, W.-K. To, On Seshadri constants of canonical bundles of compact complex hyperbolic spaces. \textit{Compositio Math.} \textbf{118} (1999), 203-215

\bibitem[Kli03]{Klingler} B. Klingler, Sur la rigidit\'e de certains groupes fundamentaux, l'arithm\'eticit\'e des r\'eseaux hyperboliques complexes, et les faux planes projectives. \textit{Invent. Math.} \textbf{153} (2003), 105-143.

\bibitem[Laz04]{Laz1} R. Lazarsfeld, Positivity in algebraic geometry. I. Classical setting: line bundles and linear series.
Ergebnisse der Mathematik und ihrer Grenzgebiete. 3. Folge. A series of Modern Surveys in Mathemathics [Results in Mathematics and Related Areas. 3rd Series. A Series of Modern Surveys in Mathematics], 48. \textit{Springer-Verlag, Berlin,} 2004.

\bibitem[MT15]{Moller-Toledo} M. M\"oller, D. Toledo, Bounded negativity of self-intersection numbers of Shimura curves in Shimura surfaces. \textit{Algebra Number Theory} \textbf{9} (2015), no. 4, 897-912.

\bibitem[PY07]{Prasad} G. Prasad, S.-K. Yeung, Fake projective planes.
\textit{Invent. Math.} \textbf{168} (2007), 321-370.

\bibitem[PY10]{Prasad2} G. Prasad, S.-K. Yeung, Addendum to ``Fake projective planes'' Invent. Math. 168, 321-370 (2007). \textit{Invent. Math.} \textbf{182} (2010), no. 1, 213-227.

\bibitem[Tol89]{Toledo} D. Toledo, Representations of surface groups in complex hyperbolic space. \textit{J. Differential Geom.} \textbf{29} (1989). no. 1, 125-133.

\bibitem[Sto15]{Stover} M.  Stover, Notes on Toledo invariant. \textit{Available on the author personal web page}.

\bibitem[Yau77]{Yau} S.-T. Yau, Calabi's conjecture and some new results in algebraic geometry. \textit{Proc. Nat. Acad. Sci. U.S.A.} \textbf{74} (1977), no. 5, 1798-1799.

\bibitem[Yeu04]{Yeung} S.-K. Yeung, Integrality and arithmeticity of co-compact lattices corresponding to certain complex two-ball quotients of Picard number one. \textit{Asian J. Math.} \textbf{8} (2004), 107-130.

\end{thebibliography}
\end{document}